\documentclass[11pt]{article}

\usepackage{graphicx}

\usepackage{hyperref}
\hypersetup{colorlinks,%
citecolor=blue,%
filecolor=blue,%
linkcolor=red,%
urlcolor=red,%
pdftex}
\usepackage{amsmath,amsfonts,amssymb}
\usepackage{amsmath}
\usepackage{multirow}
\usepackage{wrapfig}
\usepackage{epsfig}
\def\disp{\displaystyle}

\def\tto{\;{\lower 1pt \hbox{$\rightarrow$}}\kern -10pt
\hbox{\raise 2pt \hbox{$\rightarrow$}}\;}
\def\Hat{\widehat}

\def\ra{\rangle}
\def\la{\langle}

\def\B{I\!\!B}
\def\h{\hfill\square}
\def\R{I\!\!R}
\def\N{I\!\!N}
\def\ox{\bar{x}}

\def\co{\mbox{\rm co}\,}

\def\h{\hfill\square}

\def\lm{\lambda}

\def\al{\alpha}

\def \N{I\!\!N}

\newcounter{lk}

\begin{document}
\begin{center}
\vspace*{0.3in} {\bf APPLICATIONS OF CONVEX ANALYSIS TO THE\\ SMALLEST INTERSECTING BALL PROBLEM}\\[2ex]
NGUYEN MAU NAM,\footnote{Department of
Mathematics, The University of Texas--Pan American, Edinburg, TX
78539--2999, USA (email: nguyenmn@utpa.edu).} NGUYEN THAI AN \footnote{Department of
Mathematics, Hue University, 32 Leloi Hue, Vietnam (email: thaian2784@gmail.com).} and JUAN
SALINAS\footnote{Department of Mathematics, The University of
Texas--Pan American, Edinburg, TX 78539--2999, USA (email:
jsalinasn@broncs.utpa.edu).}\\[2ex]
\end{center}
\small{\bf Abstract:} \emph{The smallest enclosing circle problem}
asks for the circle of smallest radius enclosing a given set of
finite points on the plane. This problem was introduced in the 19th
century by Sylvester \cite{syl}. After more than a century, the
problem remains very active. This paper is the continuation of our
effort in shedding new light to classical geometry problems using
advanced tools of convex analysis and optimization. We propose
and study the following generalized version of the smallest
enclosing circle problem: given a finite number of nonempty closed
convex sets in a reflexive Banach space, find a ball with the
smallest radius that intersects all of the sets.

\medskip
\vspace*{0,05in} {\bf Key words.}  Convex analysis and
optimization, generalized differentiation, smallest enclosing ball problem, smallest intersecting ball problem, subgradient-type algorithms.


\newtheorem{Theorem}{Theorem}[section]
\newtheorem{Proposition}[Theorem]{Proposition}
\newtheorem{Remark}[Theorem]{Remark}
\newtheorem{Lemma}[Theorem]{Lemma}
\newtheorem{Corollary}[Theorem]{Corollary}
\newtheorem{Definition}[Theorem]{Definition}
\newtheorem{Example}[Theorem]{Example}
\renewcommand{\theequation}{\thesection.\arabic{equation}}
\normalsize

\section{Introduction and Problem Formulation}
\setcounter{equation}{0} A more general form of the smallest
enclosing circle problem is called {\em the smallest enclosing ball
problem}. Given a set $P=\{p_1, \ldots, p_n\}$, $n>1,$ on a Banach
space $X$, it is always possible to find a ball $\B(a; r)$ such that
\begin{equation*}
P\subset \B(a; r).
\end{equation*}
The smallest enclosing ball problem asks for such a ball with the
smallest radius. \vspace*{0.05in}

Consider the following optimization problem
\begin{equation}\label{OSEC}
\mbox{ minimize } f(x),\; x\in X,
\end{equation}
where the function $f$ therein is defined by
\begin{equation*}\label{SEC}
f(x)=\max\{ ||x- p_i||: i=1, \ldots, n\}.
\end{equation*}
The smallest enclosing ball can be found by solving (\ref{OSEC}).
\vspace*{0.05in}

Numerous articles have been written to study the smallest enclosing ball problem as well as its generalizations from both numerical and theoretical viewpoints. The reader are referred to \cite{chm,ljc,frank,wel} and the reference therein for recent developments as well as the history of the problem.

\begin{wrapfigure}{r}{0in}
  %\centering
  \includegraphics[width=2.9in]{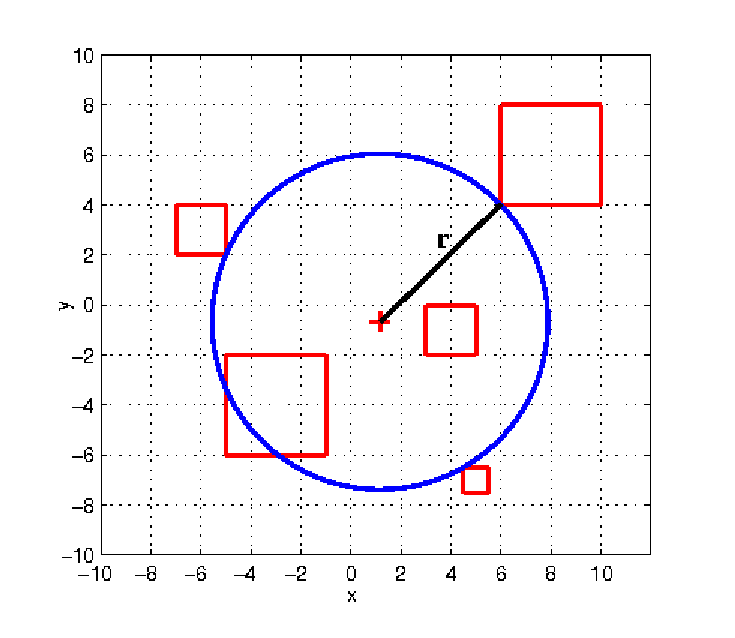}
  \caption{\scriptsize{A Smallest Intersecting Ball Problem.}}
 % \caption{Lisa Rogers}
\end{wrapfigure}

In this paper we propose and study a new problem called \emph{the
smallest intersecting ball problem} as follows: given a finite
number of nonempty closed convex sets in a Banach space, find a
ball with the smallest radius that intersects all of the sets. It is
obvious that when the sets under consideration are singleton, the
problem reduces to the smallest enclosing ball problem. This is the
continuation of our previous work from \cite{mns,mnft} in an effort
to shed new light to classical geometry problems using new tools of
nonsmooth analysis.  \vspace*{0.05in}

Let $\Omega_i, i=1, \ldots,n,$ $n>1$, be nonempty closed convex sets in
a Banach space $X$. Let $x$ be any point in $X$. Then there always
exists $r>0$ such that
\begin{equation}\label{smallest}
\B(x; r)\cap \Omega_i\neq \emptyset \mbox{ for all }i=1, \ldots,n.
\end{equation}
We are looking for a ball with the smallest radius $r>0$ (if exists)
such that property (\ref{smallest}) holds. Define the function
\begin{equation}\label{distance}
D(x)=\max\{ d(x; \Omega_i): \; i=1,\ldots,n\},
\end{equation}
where the distance function generated by a set $\Omega$ is given by
\begin{equation}\label{distance function}
d(x;\Omega)=\inf \{||x-\omega||: \omega\in \Omega\}.
\end{equation}
It is not hard to see that the function $D$ is convex since each function $d(x;\Omega_i)$, $i=1,\ldots,n$, is convex. Consider the problem of minimizing the function $D$ on $X$ below
\begin{equation}\label{SCB}
\mbox{minimize } D(x) \mbox{ subject to }x\in X.
\end{equation}
As we will see in Section 3, if $X$ is a reflexive space and $\Omega_i$,
$i=1,\ldots,n$, have no point in common, such a smallest
intersecting ball can be found by solving problem (\ref{SCB}). For
this reason, we are going to use the following standing assumptions
throughout the paper unless otherwise specified:\\[1ex]
{\em $X$ is a reflexive Banach space and $\Omega_i$, $i=1,\ldots,n,$
$n>1$, are nonempty closed convex sets in $X$ with
$$\cap_{i=1}^n\Omega_i=\emptyset.$$} Notice that if the $\Omega_i$,
$i=1,\ldots,n$, have a common point, then the smallest intersecting
ball problem has no solution (unless balls of radius $0$ are
allowed). \vspace*{0.05in}

The smallest intersecting ball problem is an example of facility
location problems.  In contrast to most of the existing facility
location problems which deal with locations of negligible sizes
(points), this new problem deals with those that involve locations
of non-negligible sizes (sets). The problem is obviously
mathematically interesting with promising applications to location
models in which it is possible to \emph{access} the entire of each
location from a point in it. The difficulty when dealing with the
smallest intersecting ball problem comes from the nonsmooth nature
of the cost function $D$ in (\ref{SCB}), especially when the norm in
$X$ is Non-Euclidean. Our approach in this paper is to study the
problem from both theoretical and numerical viewpoints using new
tools from convex analysis and optimization. The results we are
going to present in this paper provide improvements and
generalizations of many results in \cite{chm,ljc}. \vspace*{0.05in}

We organize the paper as follows. Section 2 provides necessary tools
from convex analysis and optimization for solving the smallest
intersecting ball problem. In Section 3, we study the problem from
theoretical aspects. Section 4 is devoted to developing an algorithm
of subgradient type to solve the smallest intersecting problem in
finite dimensions. The MATLAB implementations of the algorithm are
also presented.

\section{Tools of Convex Analysis}
\setcounter{equation}{0}

This section provides important constructions and results from convex analysis that will be used in the next sections.
Most of the material presented here can be found in \cite{bv,HU,r}.

Let $X$ be a normed space with the dual space $X^*$. A function $f: X\to \R$ is called convex if
\begin{equation}\label{convexity}
f\big(\lm x+(1-\lm)y\big)\le\lm f(x)+(1-\lm)f(y)\;\mbox{ for all
}\;x, y\in X\;\mbox{ and }\;\lm\in (0,1).
\end{equation}
If the inequality in (\ref{convexity}) becomes strict for $x\neq y$, the function $f$ is called strictly convex.

A subset $\Omega$ of $X$ is called convex if
\begin{equation*}
\lambda x+(1-\lambda) y\in \Omega \mbox{ for all }x, y\in \Omega \mbox{ and }\lambda\in (0,1).
\end{equation*}
It is not hard to prove that a nonempty closed subset $\Omega$ of $X$ is
convex if and only if the corresponding distance function
(\ref{distance function}) is convex. Note that the distance function
$f(x)=d(x;\Omega)$ is Lipschitz continuous on $X$ with modulus one,
i.e.,
\begin{equation*}
|f(x)-f(y)|\le||x-y||\;\mbox{ for all }\;x,y\in X.
\end{equation*}
Let $\la \cdot, \cdot\ra$ be the dual pair between $X$ and $X^*$. An element $x^*\in X^*$ is called a {\em subgradient} of a convex
function $f$ at $\ox$ if the following holds
\begin{equation*}\label{convex subdifferential}
\la x^*, x-\ox\ra\leq f(x)-f(\ox)\;\mbox{ for all }\;x\in X.
\end{equation*}
The set of all subgradients of $f$ at $\ox$ is called the \emph{subdifferential} of $f$ at $\ox$ denoted by $\partial f(\ox)$.

Convex functions and subdifferentials have several important properties as far as optimization is concerned. For instance, a convex function $f$ has a local minimum at $\ox$ if and only if it has an absolute minimum at $\ox$. Furthermore, the following generalized version of the Fermat rule holds
\begin{equation}\label{fermat}
\ox\;\mbox{ is a minimizer of }\; f\;\mbox{ if and only if
}\;0\in\partial f(\ox).
\end{equation}
It is well known that the subdifferential of the distance function
(\ref{distance function}) has a close connection to the {\em normal
cone} of the generating set $\Omega$. Recall that the normal cone of a
convex set $\Omega$ at $\ox\in \Omega$ is defined by
\begin{equation}\label{cnc}
N(\ox;\Omega)=\big\{v\in X^*:\;\la v,x-\ox\ra\le 0\;\mbox{ for all
}\;x\in\Omega\big\}.
\end{equation}
The\emph{ projection}  from a point $\ox\in X$ to a set $\Omega$ is
\begin{eqnarray}\label{pr}
\Pi(\ox;\Omega)=\big\{\omega\in\Omega\;: \;||\ox-\omega ||=d(\ox;\Omega)\big\}.
\end{eqnarray}
The following representation of subdifferential  for the distance
function (\ref{distance function}) will play an important role in
the subsequent sections of the paper. The proof of the formulas can
be found in \cite{bur}, while their various extensions are presented
in \cite{bmn10}.

\begin{Proposition}\label{subdis} Let $\Omega$ be a nonempty closed convex set of a Banach space $X$ and let $\ox\in X$. Suppose that $\Pi(\ox;\Omega)\neq \emptyset$ (which is always the case when $X$ is reflexive). Then
\begin{eqnarray*}
\partial \mbox{d}(\ox;\Omega)=\left\{\begin{array}{lr}
\partial p(\ox-\bar\omega)\cap N(\bar\omega; \Omega) &\mbox{ if
}\;\ox\notin\Omega,\\\\
N(\ox;\Omega)\cap\B^*&\mbox{ if }\;\ox\in\Omega,
\end{array}
\right.
\end{eqnarray*}
where $\B^*$ is the closed unit ball of $X^*$ and $\bar\omega$ is \emph{any element} of $\Pi(\ox;\Omega)$.

In particular, if $X$ is a Hilbert space, then $\Pi(\ox;\Omega)$ is singleton and
\begin{eqnarray*}
\partial \mbox{d}(\ox;\Omega)=\left\{\begin{array}{lr}
\big\{\dfrac{\ox-\Pi(\ox;\Omega)}{d(\ox;\Omega)}\big\} &\mbox{ if
}\;\ox\notin\Omega,\\\\
N(\ox;\Omega)\cap\B &\mbox{ if }\;\ox\in\Omega,
\end{array}
\right.
\end{eqnarray*}
\end{Proposition}
In this proposition we also observe that when $\ox\in \Omega$, one has $\Pi(\ox;\Omega)=\{\ox\}$. Since $\partial p(0)=\B^*$, the formula
\begin{equation*}
\partial d(\ox;\Omega)=\partial p(\ox-\bar\omega)\cap N(\bar\omega; \Omega)
\end{equation*}
holds for any $\ox\in X$.

Finally, we present the following well-known subdifferential rule that involves ``max" functions.
\begin{Theorem}\label{max} Let $X$ be a Banach space and let $f_i: X\to \R$, $i=1, \ldots,n,$ be continuous convex functions. Define
\begin{equation*}
f(x)=\max\{ f_i(x): i=1,\ldots, n\}.
\end{equation*}
Then
\begin{equation*}
\partial f(\ox)=\co\{ \partial f_i(\ox): i\in I(\ox)\},
\end{equation*}
where $I(\ox)=\{i=1,\ldots,n: f(\ox)=f_i(\ox)\}.$
\end{Theorem}

\section{The Smallest Intersecting Ball Problem: Theoretical \\Aspects}
\setcounter{equation}{0}

This section is devoted to theoretical analysis of the smallest intersecting ball problem.
We are able to provide improvements and generalizations of many results in \cite{chm,ljc}.
Our approach is based mostly on tools of convex analysis and optimization.

The following proposition allows us to reduce the smallest
intersecting ball problem to a nonsmooth convex optimization problem
in the reflexive space setting. For this reason, we will identify the smallest intersecting ball problem with problem (\ref{SCB}).

\begin{Proposition} Consider the minimization problem (\ref{SCB}). Then $\ox\in X$ is an
optimal solution of this problem with $r=D(\ox)$ if and only if $\B (\ox; r)$ is a smallest ball that satisfies (\ref{smallest}).
\end{Proposition}
{\bf Proof:} Suppose that $\ox$ is an optimal solution of
(\ref{SCB}) with $r=D(\ox)$. Since $\Omega_i$, $i=1,\ldots,n$, have no
point in common as in the standing assumptions, one has
\begin{equation*}
D(\ox)=\inf\{D(x): x\in X\}=r>0.
\end{equation*}
This implies
\begin{equation*}
d(\ox; \Omega_i)\leq r \mbox{ for all }i=1,\ldots, n.
\end{equation*}
Since $X$ is reflexive, there exist $\bar\omega_i\in \Omega_i, i=1,
\ldots, n$, satisfying
\begin{equation*}
||\ox-\bar\omega_i||\leq r.
\end{equation*}
It follows that $\bar\omega_i\in \B(\ox; r)\cap \Omega_i$, and hence
\begin{equation*}
\B(\ox; r)\cap \Omega_i\neq \emptyset \mbox{ for all }i=1,\ldots, n.
\end{equation*}
Suppose there exists $r'<r$ and $\bar x'\in X$ with
\begin{equation*}
\B(\bar x'; r')\cap \Omega_i\neq \emptyset \mbox{ for all } i=1, \ldots, n.
\end{equation*}
Then
\begin{equation*}
d(\bar x'; \Omega_i)\leq  r'<r \mbox{ for all } i=1,\ldots, n.
\end{equation*}
This implies $D(\bar x')\leq r'<r= D(\ox)$, which is a
contradiction. Thus $\B(\ox; r)$ is a smallest ball we are looking
for.

We are now going to justify the converse. Let us first prove that
$r=D(\ox)$. Since $\B(\ox;r)\cap \Omega_i\neq \emptyset$, one has
\begin{equation*}
d(\ox;\Omega_i)\leq r \mbox{ for all }i=1,\ldots,n.
\end{equation*}
This implies $D(\ox)\leq r$. Assume by contradiction that $D(\ox)<r$. Let $r'$ satisfy $D(\ox)<r'<r$. Then
\begin{equation*}
\B(\ox; r')\cap \Omega_i\neq \emptyset \mbox{ for all }i=1,\ldots,n.
\end{equation*}
This contradicts the minimal property of $r$. Thus $r=D(\ox)$. Let $x$ be any point in $X$ and let $r'=D(x)$. Then
\begin{equation*}
\B(x; r')\cap \Omega_i\neq \emptyset \mbox{ for all } i=1,\ldots, n.
\end{equation*}
This implies $r\leq r'$ or $D(\ox)\leq D(x')$. Therefore, $\ox$ is an optimal solution of $(\ref{SCB})$. $\h$ \vspace*{0.05in}

In what follows, we will prove that under natural assumptions on the
sets $\Omega_i$, $i=1,\ldots,n$, such a smallest intersecting ball does
exist. We are going to use the fact that on a Banach space, any
convex lower semicontinuous function is weakly lower semicontinuous.

\begin{Proposition}\label{existence} Suppose that there exists $i=1,\ldots, n,$ such that $\Omega_i$ is bounded. Then the smallest intersecting ball problem (\ref{SCB}) has a solution.
\end{Proposition}
{\bf Proof:} Without loss of generality, suppose that $\Omega_1$ is bounded. Define
\begin{equation*}
r=\inf\{ D(x): x\in X\}.
\end{equation*}
Let $(x_n)$ be a minimizing sequence for problem $(\ref{SCB})$. That means
\begin{equation*}
D(x_n)\to r \mbox{ as } n\to\infty.
\end{equation*}
Let $N\in \N$ with
\begin{equation*}
d(x_n;\Omega_1)\leq D(x_n)< r+1 \mbox{ for all } n\geq N.
\end{equation*}
Then there exists a sequence $(\omega_n)$ in $\Omega_1$ such that
\begin{equation*}
||x_n-\omega_n||< r+1 \mbox{ for all } n\geq N.
\end{equation*}
Since $(\omega_n)$ is a bounded sequence, $(x_n)$ is also bounded.
As $X$ is reflexive, there exists a subsequence $(x_{n_k})$ that
converges weakly to $\ox$. This implies
\begin{equation*}
D(\ox)\leq \liminf D(x_{n_k}) \leq r
\end{equation*}
because $D$ is weakly lower semicontinuous. Therefore, $\ox$ is a solution of problem (\ref{SCB}). $\h$ \vspace*{0.05in}

Proposition \ref{existence} implies that the smallest enclosing ball
problem (\ref{OSEC}) always has a solution because each
$\Omega_i=\{\omega_i\}$, $i=1,\ldots,n$, is obviously bounded. However,
in general, the smallest intersecting ball problem (\ref{SCB}) may
not have any solution.
\begin{Example}{\rm  Let $X=\R^2$ with the Euclidean norm. Consider $\Omega_1=\{0\}\times \R$ and $$\Omega_2=\{(x,y)\in \R^2:  y\geq \dfrac{1}{x}, x> 0\}.$$ Then the smallest intersecting ball problem (\ref{SCB}) generated by $\Omega_1$ and $\Omega_2$ does not have any solution.}
\end{Example}

In the case where the smallest intersecting ball (\ref{SCB}) has a solution, the solution may not be unique as shown in the example below.

\begin{Example}{\rm Let $X=\R^2$ with the Euclidean norm. Consider $\Omega_1=\{(x,y)\in\R^2: y\geq 1\}$ and $\Omega_1=\{(x,y)\in\R^2: y\leq -1\}$. Then any $x\in \R\times\{0\}$ is a solution of the smallest intersecting ball problem (\ref{SCB}) generated by $\Omega_1$ and $\Omega_2$.}
\end{Example}

\begin{Example} {\rm Consider $X=\R^2$ with the ``max" norm $||(x_1, x_2)||=\max\{|x_1|, |x_2|\}$. Then the
ball $\B(x; r)$ in $X$, where $x=(x_1, x_2)$ and $r>0$, is the
square
\begin{equation*}
S(x; r)=[x_1-r, x_1+r]\times [x_2-r, x_2+r].
\end{equation*}
Problem (\ref{SCB}) can be equivalently interpreted as follow: find
a smallest square $S(x;r)$ that intersects $\Omega_i$ for all
$i=1,\ldots,n$. Using different norms on $X$, we obtain different
intersecting ball problems.}
\end{Example}

\begin{Lemma}\label{sc} Let $X$ be a Hilbert space and let $\omega_i\in X$, $i=1,\ldots, n$. Then the function
\begin{equation*}
s(x)=\max\{ ||x-\omega_i||^2: i=1,\ldots, n\},\; n\geq 1,
\end{equation*}
is strictly convex.
\end{Lemma}
{\bf Proof:} We are going to prove that for $x\neq y$ and $t\in
(0,1)$, one has
\begin{equation*}
s(tx+(1-t)y)<t s(x)+ (1-t) s(y).
\end{equation*}
By induction, we only need to show that the function $p(x)=||x||^2$
is strictly convex and the function $g(x)=\max\{g_1(x), g_2(x)\}$ is
strictly convex if both $g_1$ and $g_2$ are strictly convex
functions. Indeed, for $t\in (0,1)$ and $x,y\in X$, one has
\begin{align*}
p(tx+(1-t)y)&=||tx+(1-t)y||^2\\
&=t^2||x||^2+2t(1-t)\la x, y\ra +(1-t)^2||y||^2\\
&\leq t^2||x||^2 +2t(1-t)||x||.||y|| + (1-t)^2||y||^2\\
&\leq t^2||x||^2 +t(1-t)(||x||^2+|y||^2)+(1-t)^2||y||^2\\
& = t||x||^2+(1-t)||y||^2 = tp(x)+(1-t)p(y).
\end{align*}
Notice that the equality holds if and only if $||x||=||y||$ and $\la x, y\ra =||x|| ||y||$. This implies $||x-y||^2=0$, and hence $x=y$. Therefore, $p$ is strictly convex.

Now let $t\in (0,1)$ and $x\neq y$. Then
\begin{equation*}
g_i( tx+ (1-t)y) < tg_i(x) +(1-t) g_i(y) \leq t g(x)+(1-t)g(y) \mbox{ for }i=1,2.
\end{equation*}
This implies
\begin{equation*}
g(tx+(1-t)y)<tg(x) + (1-t)g(y).
\end{equation*}
The proof is now complete. $\h$ \vspace*{0.05in}

The following proposition gives an example of a smallest intersecting ball problem which has a unique solution. We will use a natural convention that $\B(c; 0)=\{c\}$ for any $c\in X$.

\begin{Proposition}\label{uniqueness1} Let $X$ be a Hilbert space and let $r\geq 0$. Suppose that
$\Omega_i=\B(\omega_i; r)$, $i=1,\ldots,n$, are closed balls in a
Hilbert space $X$. Then the smallest intersecting ball problem
(\ref{SCB}) generated by $\Omega_i$, $i=1,\ldots,n$, has a unique
solution. Moreover, this unique solution coincides with the unique
solution of the smallest enclosing ball problem (\ref{OSEC})
generated by the centers of the balls $\{\omega_i\}$, $i=1, \ldots,
n$.
\end{Proposition}
{\bf Proof:} Let us first show that in this case the function $D$ in
(\ref{distance}) has the following representation
\begin{equation}\label{rep}
D(x) =\max \{ ||x-\omega_i||: i=1,\ldots, n\}-r.
\end{equation}
Indeed, let
\begin{equation*}
J(x)=\{i\in 1,\ldots, n: x\notin \Omega_i\}.
\end{equation*}
Since $\Omega_i$, $i=1,\ldots,n$, have no point in common by the
standing assumptions, $J(x)\neq\emptyset$. For any $i\in \{1,\ldots,
n\}\setminus J(x)$ and for any $j\in J(x)$, one has
\begin{equation*}
||x-\omega_i||\leq r \leq ||x-\omega_j||.
\end{equation*}
It follows that
\begin{align*}
D(x)&=\max \{ d(x; \Omega_i): i=1,\ldots, n\}\\
& =\max\{d(x; \Omega_i):  i\in J(x)\}\\
&=\max\{ ||x-\omega_i||: i\in J(x)\}-r\\
&=\max\{ ||x-\omega_i||: i=1,\ldots, n\}-r.
\end{align*}
Thus (\ref{rep}) has been justified. Using representation
(\ref{rep}), we see that $\ox$ is a solution of problem (\ref{SCB})
if and only if it is a solution of the minimization problem
\begin{equation}\label{square}
\mbox{ minimize } s(x)=\max\{ ||x-\omega_i||^2: i=1, \ldots, n\}, x\in X.
\end{equation}
Since $s$ is strictly convex by Lemma \ref{sc}, problem (\ref{square}) has a unique solution. Therefore, problem (\ref{SCB}) also has a unique solution. Notice that $\ox$ is a solution of problem
(\ref{square}) if and only if it is the solution of the smallest enclosing ball problem (\ref{OSEC}) generated by $\{\omega_i\}$, $i=1,\ldots,n.$ The proof is now complete. $\h$

\begin{Example}{\rm In $\R^2$ with the Euclidean norm, consider the balls $\Omega_1=\B((0,3); 3)$,
 $\Omega_2=\B((-2, 0); 1)$, and $\Omega_3=\{B((2,0); 1)\}$. Then $\ox=(0,0)$ is the
solution of problem (\ref{SCB}) generated by these balls, but this
solution is the solution of the smallest enclosing ball problem
(\ref{OSEC}) generated by the centers of the balls.}
\end{Example}

In what follows we are going to prove that in the Hilbert space
setting, the smallest intersecting ball problem (\ref{SCB})
generated by closed balls with different radii also has a unique
solution although the solution may not coincide with the solution of
the smallest enclosing ball problem generated by their centers.

\begin{Proposition} Let $X$ be a Hilbert space and let $\Omega_i=\B(\omega_i; r_i)$, $r_i\geq 0$, $i=1,\ldots,n$, be closed
balls in a Hilbert space $X$. Then the smallest intersecting ball
problem (\ref{SCB}) generated by $\Omega_i$, $i=1,\ldots, n$, has a
unique solution.
\end{Proposition}
{\bf Proof:} Similar to the proof of Proposition \ref{uniqueness1},
one has
\begin{equation*}
D(x)=\max\{||x-\omega_i||-r_i: i=1, \ldots, n\} \mbox{ for all }x\in
X.
\end{equation*}
Define
\begin{equation*}
p_i(x)=||x-\omega_i||-r_i, i=1,\ldots, n.
\end{equation*}
 Let
\begin{equation*}
\ell =-\max\{r_i: i=1,
\ldots, n\}.
\end{equation*}
Then $p_i(x)\geq \ell$ for all $i=1, \ldots, n$ and for all $x\in X$.
Consider the optimization problem
\begin{equation}\label{square1}
\mbox{minimize } h(x)=\max\{(p_i(x)-\ell)^2: i=1, \ldots, n\}.
\end{equation}
Notice that $p_i(x)-\ell\geq 0$ for all $i=1,\ldots, n$ and for all
$x\in X$. Then it is not hard to see that $\ox$ is a solution of
the smallest intersecting ball problem (\ref{SCB}) if and only if
it is also a solution of problem (\ref{square1}). Similar to the
proof of Proposition \ref{uniqueness1}, one sees that the function
$h$ in (\ref{square1}) is strictly convex and hence problem (\ref{square1}) has a unique
solution. The proof is complete. $\h$ \vspace*{0.05in}

For each $x\in X$, the set of active indices for $D$ at $x$ is defined by
\begin{equation*}
I(x)=\{ i\in \{1,\ldots,n\} : D(x)=d(x;\Omega_i)\},
\end{equation*}
and let
\begin{equation*}\label{indicator}
A_i(x)= \partial p(x-\omega_i) \cap N(\omega_i; \Omega_i).
\end{equation*}
where $\omega_i\in \Pi(x;\Omega_i)$. Notice that the definition of $A_i(x)$ does not depend on the choice of $\omega_i$ by Proposition \ref{subdis}. It is also clear from the definition that $I(x)\neq \emptyset$ for any $x\in  X$. Moreover, if $i\in I(x)$, then $D(x)>0$ because $\cap_{i=1}^n\Omega_i=\emptyset$ as in the standing assumptions, and hence
\begin{equation*}
d(x;\Omega_i)=D(x)>0.
\end{equation*}
This implies $x\notin\Omega_i$.

\begin{Proposition}\label{NS} Consider the smallest intersecting ball problem (\ref{SCB}). Then $\ox\in X$ is an optimal solution of the problem if and only if
\begin{equation*}
0\in \co \{ A_i(\ox): i\in I(\ox)\}.
\end{equation*}
\end{Proposition}
{\bf Proof:} It follows from Theorem \ref{max} that
\begin{align*}
\partial D(\ox)&=\co \{\partial d(\ox;\Omega_i): i\in I(\ox)\}\\
& =\co \{ A_i(\ox): i\in I(\ox)\}.
\end{align*}
The result then follows from the subdifferential Fermat rule (\ref{fermat}).$\h$

\begin{Corollary}\label{ConvexHull} Let $X$ be a Hilbert space. Consider the smallest intersecting ball problem (\ref{SCB}). Then $\ox$ is a solution of the problem if and only if
\begin{equation}\label{CO}
\ox\in \co \{\bar\omega_i: i\in I(\ox)\},
\end{equation}
where $\bar\omega_i= \Pi(\ox;\Omega_i)$.

In particular, if $\Omega_i=\{ a_i\}$, $i=1,\ldots,n$, then $\ox$ is the
solution of the smallest enclosing ball problem (\ref{OSEC})
generated by $a_i$, $i=1,\ldots,n$, if and only if
\begin{equation}\label{CO1}
\ox\in \co \{ a_i: i\in I(\ox)\}.
\end{equation}
\end{Corollary}
{\bf Proof:} According to Proposition \ref{NS}, the element $\ox$ is
a solution of the smallest intersecting ball problem (\ref{SCB}) if
and only if
\begin{equation*}
0\in \co \{ A_i(\ox): i\in I(\ox)\}.
\end{equation*}
For each $i\in I(\ox)$, one has $\ox\notin\Omega_i$, and hence
\begin{equation*}
A_i(\ox)=\big\{\dfrac{\ox-\bar\omega_i}{d(\ox;\Omega_i)}\big\}=\big\{\dfrac{\ox-\bar\omega_i}{D(\ox)}\big\}.
\end{equation*}
It follows that
\begin{equation*}
0\in \co \{ A_i(\ox): i\in I(\ox)\}%=\co \{ A_i(\ox): i\in I(\ox)\}
\end{equation*}
if and only if there exists $\lambda_i\geq 0$, $i \in I(\ox)$, such that $\sum_{i\in I(\ox)}\lambda_i=1$ and
\begin{equation*}
0= \sum_{i\in I(\ox)} \lambda_i \dfrac{\bar x-\bar\omega_i}{D(\ox)}.
\end{equation*}
This equation is equivalent to
\begin{equation*}
0=\sum_{i\in I(\ox)} \lambda_i(\bar x-\bar\omega_i) \mbox{ or } \ox =\sum_{i\in I(\ox)}\lambda_i\bar\omega_i,
\end{equation*}
which is equivalent to (\ref{CO}).

Notice that (\ref{CO}) is equivalent to (\ref{CO1}) when $\Omega_i=\{a_i\}$, $i=1,\ldots, n.$ The proof is now complete. $\h$ \vspace*{0.05in}

We say that the smallest ball $\B(\ox; r)$ touches a target set $\Omega_i$, $i=1,\ldots,n,$
if $\Omega_i\cap \B(\ox;r)$ is singleton.

\begin{Corollary}\label{IN} Let $X$ be a Hilbert space. Consider the smallest intersecting ball problem (\ref{SCB}). Then any smallest intersecting ball touches at least two sets among $\Omega_i$, $i=1,\ldots, n$.
\end{Corollary}
{\bf Proof:} Let us first prove that $|I(\ox)|\geq 2$ if $\ox$ is a
solution of problem (\ref{SCB}). Suppose by contradiction that
$I(\ox)=\{i_0\}$. Then by (\ref{CO}),
\begin{equation*}
\ox\in \Omega_{i_0} \mbox{ and } d(\ox; \Omega_{i_0})=D(\ox)>0.
\end{equation*}
This is a contradiction. Let us now show if $i\in I(\ox)$, then
$\B(\ox; r)$ touches $\Omega_i$, where $r=D(\ox)$. Indeed, in this case
$d(\ox; \Omega_i)=r$. If there are $u, v\in\Omega_i$ such that
\begin{equation*}
u, v\in \B(\ox; r)\cap\Omega_i,\; u\neq v.
\end{equation*}
Then
\begin{equation*}
||u-\ox||\leq r =d(\ox;\Omega_i) \mbox{ and } ||v-\ox||\leq r= d(\ox;\Omega_i).
\end{equation*}
It follows that $u, v\in \Pi(\ox;\Omega_i)$. This is a contradiction because $\Pi(\ox;\Omega_i)$ is singleton. Thus $\B(\ox; r)$ touches $\Omega_i$. The proof is now complete.
$\h$. \vspace*{0.05in}

It is obvious that $\co \{ a_i: i\in I(\ox)\}\subset \co\{a_i: i=1,
\ldots, n\}$. Thus our result in Corollary (\ref{ConvexHull}) covers
\cite[Theorem~3.6]{ljc}. It is also possible to prove that the
solution of the smallest intersecting ball problem (\ref{SCB})
generated by closed balls in a Hilbert space belongs to the convex
hull of their centers as in the next proposition.

\begin{Proposition}\label{ConvexHull1} Let $X$ be a Hilbert space.
Suppose that $\Omega_i=\B(\omega_i; r_i)$, $r_i\geq 0$, $i=1,\ldots,n$,
are closed balls in $X$. Let $\ox$ be the unique solution of problem
(\ref{SCB}). Then
\begin{equation*}
\ox\in \co\{ \omega_i: i\in I(\ox)\}.
\end{equation*}
\end{Proposition}
{\bf Proof:} Let $\ox$ be the solution of the smallest intersecting
ball problem (\ref{SCB}) generated by $\Omega_i$, $i=1,\ldots,n$. Since
$\ox\notin \Omega_i$ for all $i\in I(\ox)$, by (\ref{CO}) from Corollary
\ref{ConvexHull}, there exist $\lambda_i\geq 0$, $i\in I(\ox)$, such
that $\sum_{i\in I(\ox)}\lambda_i=1$ and
\begin{equation*}
\ox =\sum_{i\in I(\ox)} \lambda_i (\omega_i + r_i \dfrac{\ox-\omega_i}{||\ox-\omega_i||}).
\end{equation*}
For any $i\in I(\ox)$, one has $||\ox-\omega_i||=r+r_i$, where $r$ is the radius of the smallest intersecting ball. It follows that
\begin{equation*}
\ox=\sum_{i\in I(\ox)}\lambda_i\ox= \sum_{i\in I(\ox)}\lambda_i(1 -\dfrac{r_i}{r+r_i}) \omega_i + \sum_{i\in I(\ox)} \dfrac{\lambda_i r_i}{r+r_i}\ox.
\end{equation*}
This implies
\begin{equation*}
\ox=\dfrac{1}{\sum_{i\in I(\ox)} \dfrac{\lambda_i}{r+r_i}} \sum_{i\in I(\ox)} \dfrac{\lambda_i }{r+r_i} \omega_i\in\co\{\omega_i: i\in I(\ox)\}.
\end{equation*}
The proof is now complete. $\h$

\begin{Example} {\rm Let $a_i, i=1,2,3,$ be three points in $\R^2$ with the Euclidean norm and let $\ox$ be the solution of problem (\ref{OSEC}) generated by $a_i$, $i=1,2,3$. By Corollary \ref{IN}, one has $|I(\ox)|=2$ or $|I(\ox)|=3$. If $|I(\ox)|=2$, say $I(\ox)=\{2,3\}$, then $\ox\in \co\{a_2, a_3\}$ and $||\ox-a_2||=||\ox-a_3||$ by Corollary \ref{ConvexHull}. In this case
\begin{equation*}
\ox=\dfrac{a_2+a_3}{2}.
\end{equation*}
This also implies $\la a_2-a_1, a_3-a_1\ra\leq 0$. Conversely, if
$\la a_2-a_1, a_3-a_2\ra \leq 0$, then the angle of the triangle
formed by $a_1$, $a_2$, and $a_3$ at vertex $a_1$ is obtuse (we
allow the case where $a_i$, $i=1,2,3,$ are on a straight line). One
can easily see that $I(\dfrac{a_2+a_3}{2})=\{2,3\}$, and
$\dfrac{a_2+a_3}{2}$ satisfies the assumption of Corollary
\ref{ConvexHull}. Then $\ox=\dfrac{a_2+a_3}{2}$ because of the
uniqueness of the solution. In this case we have $|I(\ox)|=2$.

Thus  $|I(\ox)|=2$ if and only if one of the angles of the triangle formed by $a_1, a_2$, and $a_3$ is obtuse. In this case the solution of problem (\ref{SCB}) is the midpoint of the side opposite to the obtuse vertex.

If none of the angles of the triangle formed by $a_1, a_2, a_3$ is obtuse, then $|I(\ox)|=3$. In this case, $\ox$ is the unique point that satisfies
\begin{equation*}
\ox\in \co\{a_1, a_2, a_3\}, ||\ox-a_1||=||\ox-a_2||=||\ox-a_3||,
\end{equation*}
or $\ox$ is the center of the circumscribing circle of the triangle.}
\end{Example}
Let us now consider the solution of problem (\ref{SCB}) with the
target sets being three disjoint disks in $\R^2$.
\begin{figure}[h]
  \hfill
  \begin{minipage}[t]{.3\textwidth}
    \begin{center}
      \epsfig{file=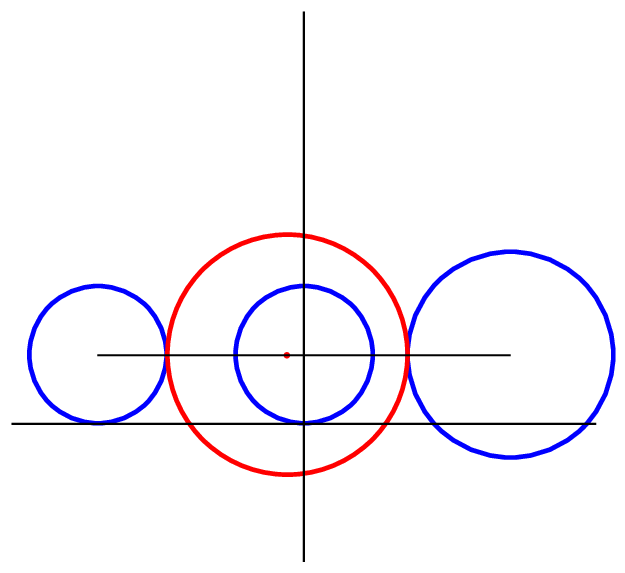, scale=0.5}
    \end{center}
  \end{minipage}
  \hfill
  \begin{minipage}[t]{.3\textwidth}
    \begin{center}
      \epsfig{file=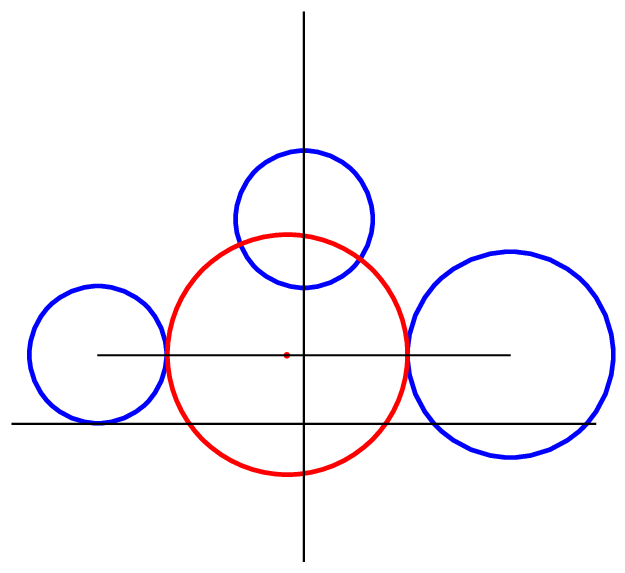, scale=0.5}
    \end{center}
  \end{minipage}
  \hfill
  \begin{minipage}[t]{.3\textwidth}
    \begin{center}
      \epsfig{file=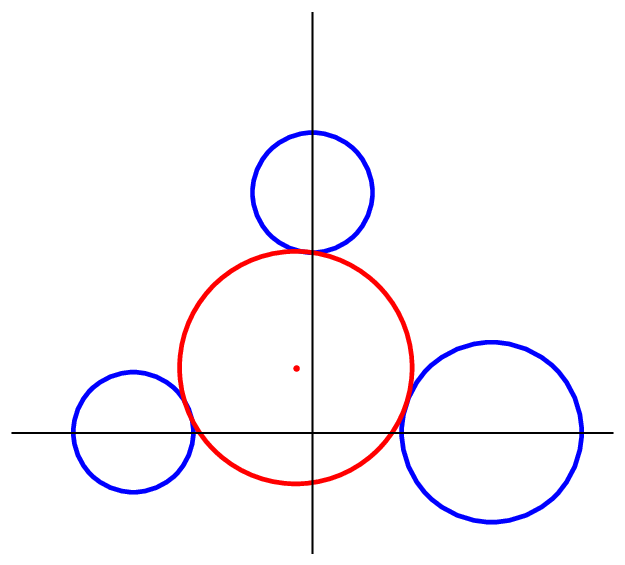, scale=0.5}
    \end{center}
  \end{minipage}
  \hfill
  \caption{Smallest Intersecting Ball Problems for Three Disks in $\R^2$}
\end{figure}

\begin{Example}{\rm Let $\Omega_i=\B(\omega_i, r_i)$, $r_i>0$, $i=1,2,3$, be disjoint disks in
$\R^2$ with the Euclidean norm. We use $bd(\Omega_i)$ to denote the boundary of $\Omega_i$, which is the
circle of center $\omega_i$ and radius $r_i$, $i=1,2,3$. Let $\ox$
be the unique solution of the problem.

Let us consider the first case where one of the line segments
connecting two of the centers intersects the other disks. For
instance, the line segment connecting $\omega_2$ and $\omega_3$
intersects $\Omega_1$. Let $u_2=\overline{\omega_2\omega_3}\cap
bd(\Omega_2)$ and $u_3=\overline{\omega_2\omega_3}\cap bd(\Omega_3)$. Let $\ox$ be the midpoint of $\overline{u_2u_3}$. Then $I(\ox)=\{2,3\}$ and we can apply  Corollary \ref{ConvexHull} to see that $\ox$ is the solution of the problem.

Now we only need to consider the case where any line segment
connecting two centers of the disks does not intersect the remaining
disk.  Let
\begin{align*}
&u_1=\overline{\omega_1\omega_2}\cap bd(\Omega_1),
v_1=\overline{\omega_1\omega_3}\cap bd(\Omega_1),\\
&u_2=\overline{\omega_2\omega_3}\cap bd(\Omega_2),
v_2=\overline{\omega_2\omega_1}\cap bd(\Omega_2),\\
&u_3=\overline{\omega_3\omega_1}\cap bd(\Omega_3),
v_3=\overline{\omega_3\omega_2}\cap bd(\Omega_3)\\
&a_1=\overline{\omega_1m_1}\cap bd(\Omega_1),
a_2=\overline{\omega_2m_2}\cap \Omega_2,
a_3=\overline{\omega_3m_3}\cap\Omega_3,
\end{align*}
where $m_1$ is the midpoint of $\overline{u_2v_3}$, $m_2$ is the
midpoint of $\overline{u_3v_1}$, and $m_3$ is the midpoint of
$\overline{u_1v_2}$. If one of the angles: $\widehat{u_2a_1v_3}$,
$\widehat{u_3a_2v_1}$, $\widehat{u_1a_3v_2}$ is greater than or
equal to $90^\circ$. For instance, if $\widehat{u_2a_1v_3}$ is
greater than or equal to $90^\circ$. Then $I(m_1)=\{2,3\}$, and
$\ox=m_1$ is the unique solution of the problem by Corollary
\ref{ConvexHull}. Now if all of the afore-mentioned angles are
acute, then $I(\ox)=3$ and the smallest disk we are looking for is
the unique disk that touches three other disks. The construction of
this disk is the celebrated problem of \emph{Apollonius}; see, e.g.,
\cite{gr}.}
\end{Example}

We are going prove that a smallest intersecting ball generated by
$n$ convex sets, $n>1$, in a $\R^m$ can be determined by at most
$m+1$ sets among them; see Figure 1 for the visualization of this
property. The proof is based on a known results for points
(see \cite[Lemma~1 (iii)]{wel}), which can be easily proved by
Corollary \ref{CO1} and the Caratheodory theorem \cite[Corollary~1,
Sec. 3.5]{IT}.

\begin{Lemma}\label{points} Let $P=\{p_1, \ldots, p_n\}$, $n>1$, be a set of finite points in $\R^m$ with the Euclidean norm
and let $\B(\ox;r)$ be the smallest enclosing ball for problem (\ref{OSEC}) generated by points in $P$. Then there exists a subset $Q\subset P$ and $2\leq |Q|\leq m+1$ such that $\B(\ox;r)$ is also the smallest enclosing ball of problem (\ref{OSEC}) generated by points in $Q$.
\end{Lemma}
{\bf Proof:} By Corollary (\ref{ConvexHull}), one has
\begin{equation*}
\ox\in \mbox{co } \{ p_i: i\in I(\ox)\}.
\end{equation*}
By the Caratheodory theorem, there exists an index set $J\subset I(\ox)$ with $|J|\leq m+1$ and
\begin{equation}\label{cv}
\ox\in \mbox{co }\{p_j: j\in J\}.
\end{equation}
It is clear that $|J|\geq 2$ because $n>1$. Let $Q=\{p_j: j\in J\}$.
Then $2\leq |Q|\leq m+1$ and $||\ox-q||=r$ for all $q\in Q$. By
converse of Corollary \ref{ConvexHull}, one has that $\B(\ox;r)$ is
the smallest enclosing ball of problem (\ref{OSEC}) generated by
points in $Q$. The proof is complete. $\h$

\begin{Proposition} Let $X=\R^m$ with the Euclidean norm. Consider problem (\ref{SCB}) in which $\Omega_i$, $i=1,\ldots,n$, are disjoint. Suppose that $\B(\ox;r)$ is a smallest intersecting ball of the problem. Then there exists an index set $J$ with $2\leq |J|\leq m+1$ such that $\B(\ox;r)$ is also a smallest intersecting ball of problem (\ref{SCB}) in which the target sets are $\Omega_j$, $j\in J$.
\end{Proposition}
{\bf Proof:} Let $\B(\ox; r)$ be a smallest intersecting ball of
problem (\ref{SCB}) with target sets $\Omega_i$, $i=1,\ldots, n$. By
Corollary \ref{ConvexHull} and Corollary \ref{IN}, one has
\begin{equation*}
\ox\in \co\{\omega_i: i\in I(\ox)\},
\end{equation*}
where $|I(\ox)|\geq 2$ and $\omega_i\in \Pi(\ox;\Omega_i)$. Again, by Corollary \ref{ConvexHull}, $\B(\ox;r)$ is the solution of the smallest enclosing ball (\ref{OSEC}) generated by $\{\omega_i: i\in I(\ox)\}$.

Applying Lemma \ref{points}, one finds a subset $J\subset I(\ox)$, $2\leq |J|\leq m+1$ such that $\B(\ox; r)$ is the solution of the smallest enclosing ball generated by $\{\omega_j: j\in J\}$. Then
\begin{equation*}
\ox\in \co\{\omega_i: i\in J\}
\end{equation*}
and $||\ox-\omega_j||=r$ for all $j\in J$. Now consider problem (\ref{SCB}) generated by $\{\Omega_j: j\in J\}$. Applying Corollary \ref{ConvexHull}, we also see that $\ox$ is a solution of this problem because
\begin{equation*}
\ox\in \co\{\omega_j: j\in J\},
\end{equation*}
where $\omega_j\in \Pi(\ox; \Omega_j)$ and $d(\ox;\Omega_j)=r$ for each $j\in J$. Moreover, $\B(\ox;r)$ is a smallest intersecting ball for the problem. The proof is now complete. $\h$ \vspace*{0.05in}

Now we are going to give a generalization of \cite[Theorem 4.4]{ljc}. Our approach, which is based on the proof of \cite[Proposition~1, Sec. 10.2]{IT}, allows us to give an estimate of the radius of the smallest intersecting ball for problem (\ref{SCB}) generated by closed balls in $\R^m$. Notice that \cite[Theorem 4.4]{ljc} holds in $\R^2$ for the classical smallest enclosing circle problem.
\begin{Theorem} Let $X=\R^m$ with the Euclidean norm.
Consider the smallest intersecting ball problem (\ref{SCB}) generated by the closed balls $\Omega_i=\B(\omega_i; r_i)$, $i=1,\ldots,n$. Let
\begin{align*}
&r_{min}=\min\{r_i: i=1, \ldots, n\}, r_{\max}:=\max\{r_i: i=1, \ldots, n\},\\
&\ell=\min\{m+1, n\}, P=\{\omega_i: i=1, \ldots, n\},
\end{align*}
and let $B_*=\B(\ox; r)$ be the smallest enclosing ball. Then
\begin{equation}\label{estimate}
\dfrac{1}{2}\mbox{diam }P-r_{max}\leq r\leq \sqrt{\dfrac{\ell-1}{2\ell}}\mbox{diam }P-r_{min}.
\end{equation}
In particular,
\begin{equation*}\label{estimate1}
\mbox{diam }P-\mbox{max }\{\mbox{diam} \Omega_i: i=1, \ldots, n\}\leq \mbox{diam } B_*\leq \sqrt{\dfrac{2(\ell-1)}{\ell}}\mbox{diam }P-\mbox{min }\{\mbox{diam} \Omega_i: i=1, \ldots, n\}.
\end{equation*}
\end{Theorem}
{\bf Proof:} For any $i,j=1,\ldots,n$, one has
\begin{equation*}
||\omega_i-\omega_j||\leq ||\ox-\omega_i||+||\ox-\omega_j||\leq 2r +
r_i+r_j\leq 2r +2 r_{max}.
\end{equation*}
Thus
\begin{equation*}
\mbox{diam }P\leq 2r+2r_{max}
\end{equation*}
and the first inequality in (\ref{estimate}) holds true.

Let us prove the second inequality of (\ref{estimate}). By the
Caratheodory theorem \cite[Corollary~1, Sec. 3.5]{IT} and the proof
of Proposition \ref{ConvexHull1}, there exist $k\leq \min\{m+1,
I(\ox)\}\leq \min\{m+1, n\}$ and $\lambda_i\geq 0$, $i=1, \ldots,
k$, $\sum_{i=1}^k\lambda_i=1$ (we reorder the indices if necessary)
such that
\begin{equation*}
\ox=\sum_{i=1}^k \mu_i \omega_i, \mbox{ where } \mu_i=\dfrac{1}{\mu}\dfrac{\lambda_i}{r_i+r}, \mu=\sum_{i=1}^k \dfrac{\lambda_i}{r+r_i}
\end{equation*}
and
\begin{equation*}
||\ox-\omega_i||=r+r_i, i=1, \ldots, k.
\end{equation*}
One has
\begin{align*}
||\omega_i-\omega_j||^2&=||\omega_i-\ox||^2 +||\omega_j-\ox||^2 -2 \la \omega_i-\ox, \omega_j-\ox\ra\\
&=(r+r_i)^2 +(r+r_j)^2-2\la \omega_i-\ox, \omega_j-\ox\ra.
\end{align*}
Thus
\begin{align*}
\sum_{i=1}^k\mu_i ||\omega_i-\omega_j||^2 &= (r+r_j)^2 +\sum_{i=1}^k \mu_i (r+r_i)^2 -2 \la \sum_{i=1}^k \mu_i\omega_i-\ox, \omega_j-\ox\ra\\
&=(r+r_j)^2 +\sum_{i=1}^k \mu_i (r+r_i)^2.
\end{align*}
It follows that
\begin{equation*}
\sum_{i, j=1}^k \mu_i\mu_j ||\omega_i-\omega_j||^2=2\sum_{i=1}^k \mu_i(r+r_i)^2.
\end{equation*}
We also have
\begin{equation*}
\sum_{i, j=1, i\neq j}^k \mu_i\mu_j=(\sum_{i=1}^k\mu_i)^2-\sum_{i=1}^k\mu_i^2\leq 1-\dfrac{1}{k}=\dfrac{k-1}{k}\leq \dfrac{\ell-1}{\ell}.
\end{equation*}
This implies
\begin{equation*}
\sum_{i, j=1}^k \mu_i\mu_j ||\omega_i-\omega_j||^2=2\sum_{i=1}^k \mu_i(r+r_i)^2\leq (\mbox{diam }P)^2\dfrac{\ell-1}{\ell}.
\end{equation*}
Since
\begin{equation*}
\mu_i(r+r_i)^2= \dfrac{1}{\mu}\dfrac{\lambda_i}{r_i+r} (r+r_i)^2=\dfrac{\lambda_i}{\mu}(r+r_i),
\end{equation*}
one has
\begin{align*}
\sum_{i=1}^k \dfrac{\lambda_i}{\mu}(r+r_i)\leq \dfrac{1}{2}\dfrac{\ell-1}{\ell}(\mbox{diam }P)^2
\end{align*}
Using the formula for $\mu$, we arrive at
\begin{equation*}
(r+r_{min})^2 \leq \dfrac{1}{2}\dfrac{\ell-1}{\ell}(\mbox{diam }P)^2.
\end{equation*}
This implies (\ref{estimate}). The second estimates follow from (\ref{estimate}). The proof is complete. $\h$

\section{Subgradient Algorithm and Its Implementation}
\setcounter{equation}{0}

In this section let $(X, ||\cdot||)$ be a normed space where
$X=\R^m$ and let $p(x)=||x||$ be the norm function on $X$. We are
going to present and justify an algorithm of  subgradient type to
solve problem (\ref{SCB}) numerically and illustrate its
implementations using MATLAB.

\begin{Theorem}\label{subgradient method2} Let $\Omega_i$,
$i=1,\ldots,n$, be nonempty closed convex subsets of $X$ such
that at least one of them is bounded. Picking a sequence
$\{\al_k\}$ of positive numbers and a starting point $x_1\in
X$, consider the iterative algorithm:
\begin{equation}\label{al}
x_{k+1}=x_k-\alpha_k x^*_k,\;
k\in\N.
\end{equation}
Let the vectors $x^*_k$ in (\ref{al}) be given by
\begin{eqnarray}\label{a1}
x^*_k\in \partial p(x_k-\omega_k)\cap N(\omega_k; \Omega_i),
\end{eqnarray}
where $\omega_k \in \Pi(x_k;\Omega_i)$ and $i$ is any index chosen from
the following index set
\begin{equation*}
I(x_k)=\{ i=1, \ldots, n: D(x_k)=d(x_k;\Omega_i)\}.
\end{equation*}
Assume that the given sequence $\{\al_k\}$ in \eqref{al} satisfies
the conditions
\begin{equation}\label{a2}
\sum_{k=1}^\infty\alpha_k=\infty\;\mbox{ and
}\;\sum_{k=1}^\infty\alpha_k^2<\infty.
\end{equation}
Then the iterative sequence $\{x_k\}$ in \eqref{a1} converges to
an optimal solution of the smallest intersecting ball problem
\eqref{SCB} and the value sequence
\begin{equation}\label{Vk}
V_k=\min\big\{D(x_j): j=1,\ldots,k\big\}
\end{equation}
converges to the optimal value $\Hat V$ in this problem.

Furthermore, we have the estimate
\begin{align*}
V_k-\Hat
V\le\dfrac{{\rm{d_2}}(x_1;S)^2+\ell^2}{2\sum_{i=1}^k\alpha_k},
\end{align*}
where $\ell^2=\sum_{k=1}^\infty\alpha_k^2$, and $d_2(x_1;S)$ denotes the distance generated by the Euclidean norm from $x_1$ to the solution set $S$ of the problem.
\end{Theorem}
{\bf Proof:} By Proposition \ref{existence} the smallest
intersecting ball problem under consideration has a solution.
Observe that the function $D$ in (\ref{SCB}) satisfies a Lipschitz
condition with Lipschitz constant $\kappa=1$. We have
\begin{align*}
\partial D(x_k)&=\co\{ \partial d(x_k; \Omega_i): i \in I(x_k)\}\\
&= \co \{ \partial p(x_k-\omega_k)\cap N(\omega_k; \Omega_i): i \in I(x_k)\},
\end{align*}
where $\omega_k\in \Pi(x_k; \Omega_i)$. Notice that under the standing
assumptions, $x_k\notin\Omega_i$ for $i\in I(x_k)$. It follows that for
any $i\in I(x_k)$ one has
\begin{equation*}
\partial p(x_k-\omega_k)\cap N(\omega_k; \Omega_i)\subset \partial D(x_k).
\end{equation*}
We also have that $\partial d(x_k; \Omega_i)=\partial p(x_k-\omega_k)\cap
N(\omega_k; \Omega_i)$ is nonempty. Since all norms in $X$ are equivalent, it
suffices to show that
\begin{equation*}
||x_k-\bar x||_2\to 0 \mbox{ and }V_k\to \Hat V,
\end{equation*}
where $||\cdot||_2$ is the Euclidean norm in $X$ and $\bar x$ is a solution of problem (\ref{SCB}).
However, these follow directly from the well-known results on the subgradient method
for convex functions in the so-called ``square summable but not
summable case"; see, e.g., \cite{boyd,bert}. $\h$ \vspace*{0.05in}

One important features of the subgradient method is that the subgradient $x^*_k$ for each $k$ is not uniquely defined. This also reflects in the following direct consequence of Theorem \ref{subgradient method2}.

\begin{Corollary} Let $X=\R^m$ with the Euclidean norm. Consider the
smallest intersecting ball problem (\ref{SCB}). For each
$k\in \N$, the subgradient $x^*_k$ in Theorem \ref{subgradient method2} is computed by
\begin{equation*}\label{indicator}
x^*_k= \dfrac{x_k-\omega_k}{||x_k-\omega_k||},
\end{equation*}
where $\omega_k= \Pi(x_k; \Omega_i)$ and $i$ is an index chosen from $I(x_k)$.

In particular, if $\Omega_i=\B(c_i; r_i)$, $i=1,\ldots,n$, are closed balls in $X$. Then the subgradient $x^*_k$ has
the following explicit representation
\begin{equation*}\label{indicator}
x^*_k= \dfrac{x_k-c_i}{||x_k-c_i||} \mbox{ for an index } i\in I(x_k).
\end{equation*}
\end{Corollary}
\begin{figure}[h]
\vspace{-25pt}
\begin{minipage}{2in}
   \vspace{15pt}
   \hspace{15pt}\includegraphics[width=4.2in]{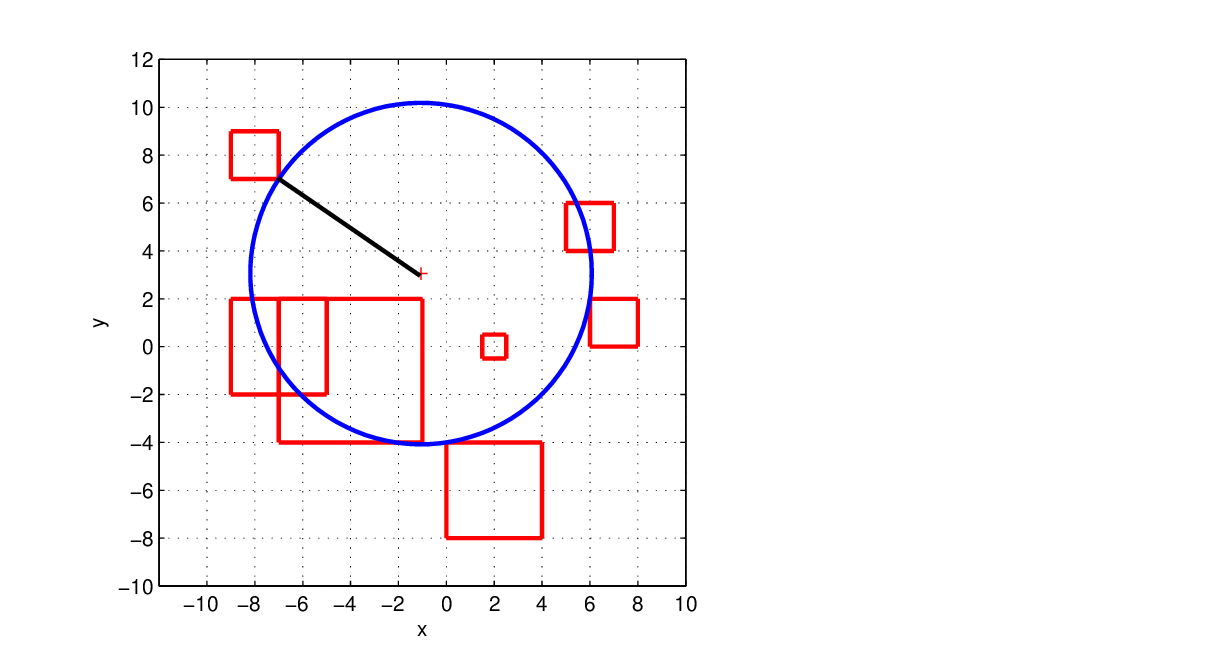}\\
\end{minipage}
~\hfill~
\begin{minipage}[t]{0.54\textwidth}
\begin{tabular}{|c|c|c|}
\hline
\multicolumn{3}{|c|}{MATLAB RESULT} \\
\hline
$k$ & $x_k$ & $V_k$ \\
\hline
1         & (2,2)              & 10.29563 \\
1,000     & (-1.05477,3.05560) & 7.13408 \\
10,000    & (-1.05562,3.05553) & 7.13408 \\
100,000   & (-1.05555,3.05555) & 7.13408 \\
400,000   & (-1.05556,3.05556) & 7.13408 \\
800,000   & (-1.05556,3.05556) & 7.13408 \\
1,000,000 & (-1.05556,3.05556) & 7.13408 \\
\hline
\end{tabular}
\end{minipage}
\vspace{-30pt} \caption{A Smallest Intersecting Ball with Euclidean
Norm to Square Targets.}
\end{figure}

\begin{Example}\label{ES}{\rm Consider $X=\R^2$ with the Euclidean norm. The target sets are the squares $\Omega_i=S(\omega_i; r_i)$, $\omega_i=(\omega_{1i}, \omega_{2i}),$ $i=1, \ldots, n$, where
\begin{equation*}
S(\omega_i; r_i)=[\omega_{1i}-r_i, \omega_{1i}+r_i]\times [\omega_{2i}-r_i, \omega_{2i}+r_i].
\end{equation*}
Let the vertices of the $i$th square be denoted by
$v_{1i}=(\omega_{1i}+r_i,\omega_{2i}+r_i),v_{2i}=(\omega_{1i}-r_i,\omega_{2i}+r_i),v_{3i}=(\omega_{1i}-r_i,\omega_{2i}-r_i),v_{4i}=(\omega_{1i}+r_i,\omega_{2i}-r_i)$
and let $x_k=(x_{1k}, x_{2k})$. Fix an index $i\in I(x_k)$. Then the vectors $x^*_k$ in Theorem
\ref{subgradient method2} are given by
{\small
\begin{equation*}
x^*_k=\left\{\begin{array}{ll}
\disp\frac{x_k-v_{1i}}{\|x_k-v_{1i}\|} &\mbox{if }\;x_{1k}-\omega_{1i}>r_i \mbox{ and }\;x_{2k}-\omega_{2i}> r_i,\\\\
\disp\frac{x_k-v_{2i}}{\|x_k-v_{2i}\|} &\mbox{if }\;x_{1k}-\omega_{1i}<-r_i \mbox{ and }\;x_{2k}-\omega_{2i}> r_i,\\\\
\disp\frac{x_k-v_{3i}}{\|x_k-v_{3i}\|} &\mbox{if }\;x_{1k}-\omega_{1i}<-r_i \mbox{ and }\;x_{2k}-\omega_{2i}< -r_i,\\\\
\disp\frac{x_k-v_{4i}}{\|x_k-v_{4i}\|} &\mbox{if }\;x_{1k}-\omega_{1i}>r_i \mbox{ and }\;x_{2k}-\omega_{2i}< -r_i,\\\\
(0,1) &\mbox{if }\;|x_{1k}-\omega_{1i}|\le r_i \mbox{ and }\;x_{2k}-\omega_{2i}> r_i,\\\\
(0,-1) &\mbox{if }\;|x_{1k}-\omega_{1i}|\le r_i \mbox{ and }\;x_{2k}-\omega_{2i}< -r_i,\\\\
(1,0) &\mbox{if }\; x_{1k}-\omega_{1i}> r_i \mbox{ and }\;|x_{2k}-\omega_{2i}|\le r_i,\\\\
(-1,0) &\mbox{if }\; x_{1k}-\omega_{1i}< -r_i \mbox{ and }\;|x_{2k}-\omega_{2i}|\le r_i.\\\\
\end{array}\right.
\end{equation*}}
It is also not hard to determine the index set $I(x_k)$ sequence and $V_k$ for each $k\in\N$. Thus the algorithm is explicit.

Consider the target sets $\Omega_i, i=1,\ldots, 7$, to be the squares with centers
$(-8,8)$, $(-7,0)$, $(-4,-1)$, $(2,0)$, $(2,-6)$, $(7,1)$, and $(6,5)$ and
the radii $r_i=\{1,2,3,0.5,2,1,1\}$ for $i=1,\ldots, 7$, respectively. A MATLAB program is performed for the sequence
$\al_k=1/k$ satisfying \eqref{a2} and the starting point $x_1=$
(2,2); see Figure 3.

Observe that the numerical results computed up to five decimal
places yield an optimal solution $\ox\approx(-1.05556,3.05556)$ and
the optimal value $\Hat V\approx7.13408$.}
\end{Example}

When working with a norm in $X$ that is different from the Euclidean
norm, it may be difficult to find the distance functions, the
projections to sets, as well as the subdifferential of the norm. The
following remark allows us to have an intuitive way to find a
subgradient $x^*_k$, $k\in\N$, in Theorem (\ref{subgradient
method2}) in the case $X=\R^2$ with the ``sum" norm.

\begin{Remark}\label{sum norm}{\rm Let $X=\R^2$ with the ``sum" norm $p(x)=|x_1|+|x_2|$, $x=(x_1,x_2)$. The ball
$\B(\ox; t)$, $\ox=(\ox_1,\ox_2)$, $t>0$, is the following diamond
shape
\begin{equation*}
\B(\ox; t)=\{(x_1, x_2)\in X: |x_1-\ox_1|+|x_2-\ox_2|\leq t\}.
\end{equation*}
The distance from $\ox$ to a nonempty closed set $\Omega$ and the corresponding projection are given by
\begin{equation}\label{minimal time}
d(\ox;\Omega)=\min\{t\geq 0: \B(\ox; t)\cap \Omega\neq \emptyset\}
\end{equation}
and
\begin{equation*}
\Pi(\ox;\Omega)=\B(\ox; t)\cap \Omega, \mbox{ where }t= d(\ox;\Omega).
\end{equation*}
Moreover, the subdifferential $\partial p(\ox)$, $\ox\in X$, has the following explicit representation
{\small\begin{eqnarray*}
\partial p(\ox_1,\ox_2)=\left\{\begin{array}{ll}
[-1,1]\times [-1,1],&\mbox{if }\;
(\ox_1,\ox_2)=(0,0),\\\\
\disp [-1,1]\times \{1\},&\mbox{if }\;\ox_1=0, \ox_2>0,\\\\
\disp [-1,1]\times \{-1\},&\mbox{if }\;\ox_1=0, \ox_2<0,\\\\
\disp \{1\}\times [-1,1],&\mbox{if }\;\ox_1>0, \ox_2=0,\\\\
\disp \{-1\}\times [-1,1],&\mbox{if }\;\ox_1<0, \ox_2=0,\\\\
\disp \{1\}\times \{1\},&\mbox{if }\;\ox_1>0, \ox_2>0,\\\\
\disp \{1\}\times \{-1\},&\mbox{if }\;\ox_1>0, \ox_2<0,\\\\
\disp \{-1\}\times \{1\},&\mbox{if }\;\ox_1<0, \ox_2>0,\\\\
\disp \{-1\}\times \{-1\},&\mbox{if }\;\ox_1<0, \ox_2<0.\\\\
\end{array}\right.
\end{eqnarray*}}}
\end{Remark}
By considering the ``sum" norm in $X$, we are able to introduce a new smallest intersecting ball problem in which a ``ball" is a diamond shape. The algorithm is going to be implemented in the following example.

\begin{Example}{\rm
Let us consider an example when $X=\R^2$ with the ``sum" norm. Let
$\Omega_i$ be the squares $S(\omega_i; r_i)$, $i=1,\ldots,n$, given in Example \ref{ES}. Notice that a
ball $\B(\ox; r)$, $\ox=(\bar x_1, \bar x_2)$, in $X$ is the diamond shape
\begin{equation*}
\B(\ox; r)=\{ (x_1, x_2)\in \R^2: |x_1-\ox_1|+|x_2-\ox_2|\leq r\}.
\end{equation*}
Therefore, the smallest intersecting ball problem
(\ref{SCB}) can be interpreted as follows: find a diamond shape in
$\R^2$ that intersects all $n$ given squares.
Using the same notation for the vertices of the target set $\Omega_i$ as
in Example \ref{ES}, one can see that the vectors $x^*_k$
in Theorem \ref{subgradient method2} are given by {\small
\begin{eqnarray*}
x^*_k=\left\{\begin{array}{ll} (1,1),&\mbox{ if
}\;x_{1k}-\omega_{1i}>r_i\;
\mbox{ and }\;x_{2k}-\omega_{2i}>r_i,\\\\
(-1,1),&\mbox{ if
}\;x_{1k}-\omega_{1i}<-r_i\;\mbox{ and }\;x_{2k}-\omega_{2i}>r_i,\\\\
(-1,-1),&\mbox{ if }\;x_{1k}-\omega_{1i}<-r_i\;
\mbox{ and }\;x_{2k}-\omega_{2i}<-r_i,\\\\
(1,-1),&\mbox{ if }\;x_{1k}-\omega_{1i}>r_i\;
\mbox{ and }\;x_{2k}-\omega_{2i}<-r_i,\\\\
(0,1), & \mbox{ if }\;|x_{1k}-\omega_{1i}|\le r_i\;\mbox{ and
}\;x_{2k}-\omega_{2i}>r_i,\\\\
(0,-1), & \mbox{ if }\;|x_{1k}-\omega_{1i}|\le r_i\;\mbox{ and }
\;x_{2k}-\omega_{2i}<-r_i,\\\\
(1,0), & \mbox{ if } \;x_{1k}-\omega_{1i}>r_i\;\mbox{ and
}\;|x_{2k}-\omega_{2i}|\le
r_i,\\\\
(-1,0), & \mbox{ if }\;x_{1k}-\omega_{1i}<-r_i\;\mbox{ and
}\;|x_{2k}-\omega_{2i}|\le r_i.
\end{array}\right.
\end{eqnarray*}}
Consider the target sets $\Omega_i, i=1,\ldots, 7$, to be the squares with centers
$(-5,3)$, $(-3,0)$, $(-2,-3)$, $(0,-8)$, $(4,-3)$, $(3,0)$, and $(5,4)$ and
the radii $r_i=1$ for $i=1,\ldots, 6$, respectively. A MATLAB program is performed for the sequence
$\al_k=1/k$ satisfying \eqref{a2} and the starting point $x_1=$
(2,0); see Figure 4.
\begin{figure}[h]
\vspace{-25pt}
\begin{minipage}{2in}
   \vspace{15pt}
   \hspace{15pt}\includegraphics[width=4.05in]{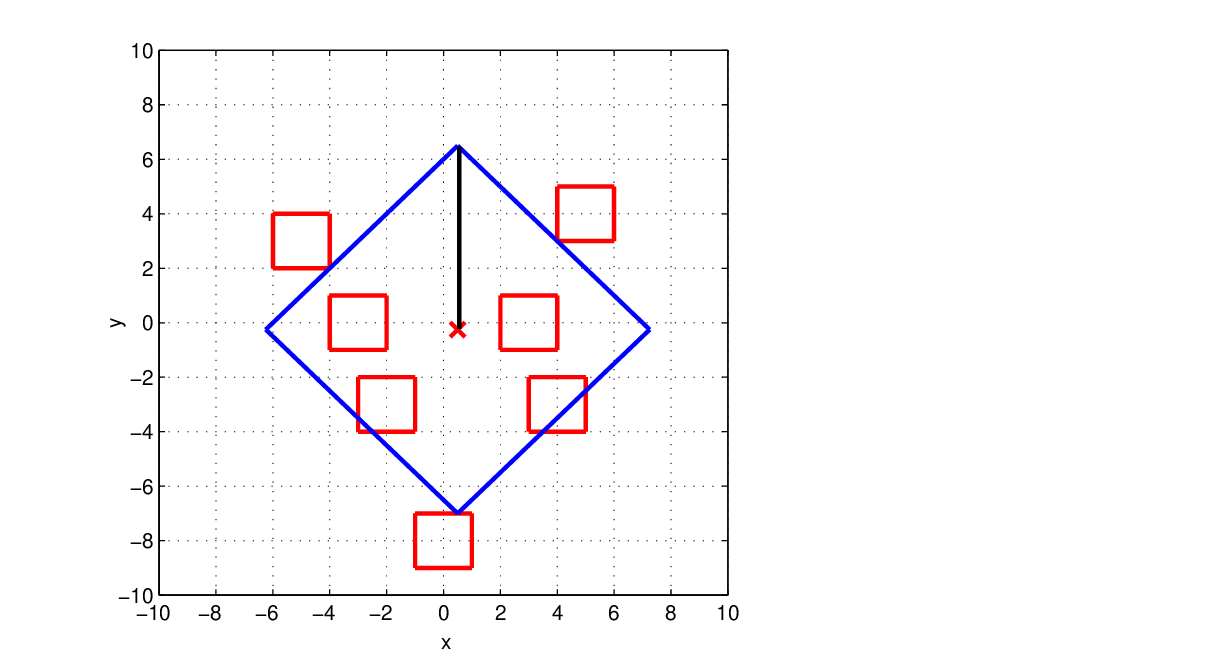}\\
\end{minipage}
~\hfill~
\begin{minipage}[t]{0.53\textwidth}
\begin{tabular}{|c|c|c|}
\hline
\multicolumn{3}{|c|}{MATLAB RESULT} \\
\hline
$k$ & $x_k$ & $V_k$ \\
\hline
1       & (2,0)              & 8 \\
10      & (0.67500,-0.17897) & 6.85397 \\
100     & (0.49999,-0.25853) & 6.75040 \\
1,000   & (0.50000,-0.25049) & 6.75040 \\
10,000  & (0.50000,-0.25004) & 6.75004 \\
100,000 & (0.50000,-0.25000) & 6.75000 \\
200,000 & (0.50000,-0.25000) & 6.75000 \\
\hline
\end{tabular}
\end{minipage}
\vspace{-30pt}
\caption{A Smallest Intersecting Ball with Sum Norm to Square Targets.}
\end{figure}

Observe that the numerical results computed up to five decimal
places yield optimal solution $\ox\approx(0.50000,-0.25000)$ and the
optimal value $\Hat V\approx6.75000$.}
\end{Example}

Similar observations for $X=\R^2$ with the ``max" norm can be easily seen:

\begin{Remark}{\rm Let $X=\R^2$ with the ``max" norm $p(x)=\max\{|x_1|,|x_2|\}$, $x=(x_1,x_2)$. The ball $\B(\ox; t)$,
$\ox=(\ox_1,\ox_2)$, $t>0$, is the following square
\begin{equation*}
\B(\ox; t)=[\ox_1-r, \ox_1+r]\times [\ox_2-r, \ox_2+r].
\end{equation*}
The distance from $\ox$ to a nonempty closed set $\Omega$ and the corresponding projection are given by
\begin{equation*}\label{minimal time}
d(\ox;\Omega)=\min\{t\geq 0: \B(\ox; t)\cap \Omega\neq \emptyset\}
\end{equation*}
and
\begin{equation*}
\Pi(\ox;\Omega)=\B(\ox; t)\cap \Omega, \mbox{ where }t= d(\ox;\Omega).
\end{equation*}
Moreover, the subdifferential $\partial p(\ox)$, $\ox\in X$, has the following explicit representation
{\small\begin{eqnarray*}
\partial p(\ox_1,\ox_2)=\left\{\begin{array}{ll}
\big\{(v_1,v_2)\in\R^2\;\big|\;|v_1|+|v_2|\le 1\big\}&\mbox{if }\;
(\ox_1,\ox_2)=(0,0),\\\\
\{(0,1)\}&\mbox{if }\;|\ox_1|<\ox_2,\\\\
\{(0,-1)\}&\mbox{if }\;\ox_2<-|\ox_1|,\\\\
\{(1,0)\}&\mbox{if }\;x_1>|\ox_2|,\\\\
\{(-1,0)\}&\mbox{if }\;\ox_1<-|\ox_2|,\\\\
\big\{(v_1,v_2)\in\R^2\;\big|\;|v_1|+|v_2|=1,\;v_1\ge 0,\;
v_2\ge 0\big\}&\mbox{if }\;\ox_1=\ox_2>0,\\\\
\big\{(v_1,v_2)\in\R^2\;\big|\;|v_1|+|v_2|=1,\;v_1\ge 0,\;v_2\le
0\big\}&\mbox{if }\;\ox_1=-\ox_2>0,\\\\
\big\{(v_1,v_2)\in\R^2\;\big|\;|v_1|+|v_2|=1,\;v_1\le 0,\;v_2\le
0\big\}&\mbox{if }\;\ox_1=\ox_2<0,\\\\
\big\{(v_1,v_2)\in\R^2\;\big|\;|v_1|+|v_2|=1,\;v_1\le 0,\;v_2\ge
0\big\}&\mbox{if }\;\ox_1=-\ox_2<0.
\end{array}\right.
\end{eqnarray*}}}

\end{Remark}

\begin{Example}{\rm
Let us consider an example when $X=\R^2$ with the ``max" norm. Let
$\Omega_i$ be the squares given in Example \ref{ES}. Notice that a
ball $\B(\ox; r)$, $\ox=(\bar x_1, \bar x_2)$, in $X$ is the square
\begin{equation*}
\B(\ox; r)=[\ox_1-r, \ox_1+r]\times [\ox_2-r, \ox_2+r].
\end{equation*}
Therefore, the smallest intersecting ball problem
(\ref{SCB}) can be interpreted as follows: find a smallest square in
$\R^2$ that intersects all $n$ given squares. Using
the same notation for the vertices of the target set $\Omega_i$ as Example \ref{ES}, one can see that the vectors $x^*_k$ in
Theorem \ref{subgradient method2} are given by {\small
\begin{eqnarray*}
x^*_k=\left\{\begin{array}{ll} (1,0),&\mbox{if }\;|x_{2k}-\omega_{2i}|\le
x_{1k}-\omega_{1i}\;\mbox{ and }\;x_{1k}>\omega_{1i}+r_i,\\\\
(-1,0),&\mbox{if }\;|x_{2k}-\omega_{2i}|\le \omega_{1i}-x_{1k}\;\mbox{ and }\;
x_{1k}<\omega_{1i}-r_i,\\\\
(0,1),&\mbox{if }\;|x_{1k}-\omega_{1i}|\le x_{2k}-\omega_{2i}\;\mbox{ and }\;
x_{2k}>\omega_{2i}+r_i,\\\\
(0,-1),&\mbox{if }\;|x_{1k}-\omega_{1i}|\le
\omega_{2i}-x_{2k}\;\mbox{ and }\; x_{2k}<\omega_{2i}-r_i,
\end{array}\right.
\end{eqnarray*}}
where $i\in I(x_k)$. The sequence $(V_k)$ is determined based on
$D(x_k)$. Fix any $i\in I(x_k)$. Then {\small\begin{eqnarray*}
D(x_{1k},x_{2k})=\left\{\begin{array}{ll}
x_{1k}-(\omega_{1i}+r_i),&\mbox{if }\;|x_{2k}-\omega_{2i}|\le x_{1i}-\omega_{1i},\;x_{1k}>\omega_{1i}+r_i,\\\\
(\omega_{1i}-r_i)-x_{1k},&\mbox{if }\;|x_{2k}-\omega_{2i}|\le \omega_{1i}-x_{1k},\;x_{1k}<\omega_{1i}-r_i,\\\\
x_{2k}-(\omega_{2i}+r_i),&\mbox{if }\;|x_{1k}-\omega_{1i}|\le x_{2k}-\omega_{2i},\;x_{2k}>\omega_{2i}+r_i,\\\\
(\omega_{2i}-r_i)-x_{2k},&\mbox{if }\;|x_{1k}-\omega_{1i}|\le
\omega_{2i}-x_{2k},\;x_{2k}< \omega{2i}-r_i.
\end{array}
\right.
\end{eqnarray*}}
Consider the target sets $\Omega_i, i=1,\ldots, 6$, to be the squares with centers
$(-5,7)$, $(-2,0)$, $(2,-5)$, $(7,-2)$, $(3,2)$, and $(7,8)$ and
the radii $r_i=\{1,1,0.5,1,2,0.5\}$ for $i=1,\ldots, 6$, respectively. A MATLAB program is performed for the sequence
$\al_k=1/k$ satisfying \eqref{a2} and the starting point $x_1=$
(-2,3); see Figure 5.
\begin{figure}[h]
\vspace{-25pt}
\begin{minipage}{2in}
   \vspace{15pt}
   \hspace{15pt}\includegraphics[width=4.05in]{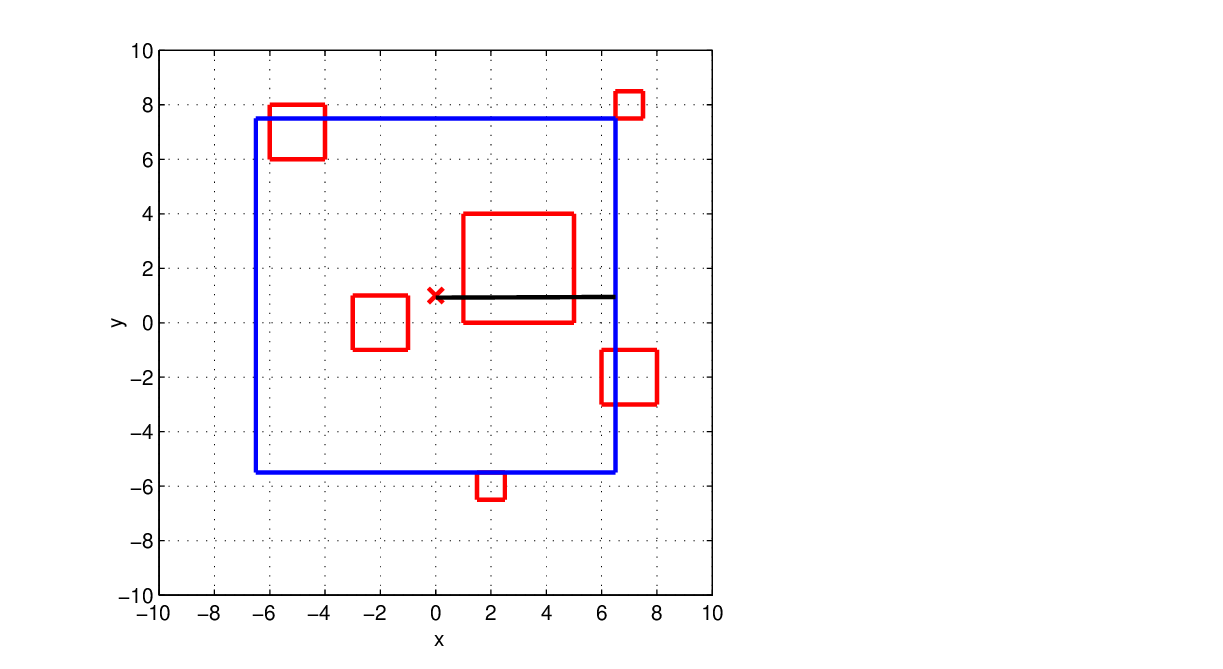}\\
\end{minipage}
~\hfill~
\begin{minipage}[t]{0.54\textwidth}
\begin{tabular}{|c|c|c|}
\hline
\multicolumn{3}{|c|}{MATLAB RESULT} \\
\hline
$k$ & $x_k$ & $V_k$ \\
\hline
1       & (-2,3)             & 8.50000 \\
10      & (-0.52500,1.54603) & 7.04603 \\
100     & (0.02973,1.00990)  & 6.50000 \\
1,000   & (0.02973,1.00100)  & 6.50000 \\
10,000  & (0.02973,1.00010)  & 6.50000 \\
100,000 & (0.02973,1.00000)  & 6.50000 \\
200,000 & (0.02973,1.00000)  & 6.50000 \\
\hline
\end{tabular}
\end{minipage}
\vspace{-30pt}
\caption{A Smallest Intersecting Ball with Max Norm to Square Targets.}
\end{figure}
Observe that the numerical results computed up to five decimal
places yield an optimal solution $\ox\approx(0.02973,1.00000)$ and the
optimal value $\Hat V\approx6.50000$.}
\end{Example}

The advantage of the algorithm comes from the fact that we are able
to deal with the smallest intersecting ball problem generated by
target sets of different types and different norms. Although a
faster subgradient algorithm may be applied to this problem, we have
chosen the simplest one for demonstrations.

\end{document}